\documentclass[a4paper,12pt]{article}

\usepackage{amsmath}
\usepackage{amsfonts}
\usepackage{amsthm}
\usepackage{amssymb}

\newtheorem{thm}{Theorem}[section]
\newtheorem{lem}[thm]{Lemma}
\newtheorem{prop}[thm]{Proposition}
\newtheorem*{mainresult}{Theorem \ref{main_theorem}}

\theoremstyle{definition}
\newtheorem{eg}[thm]{Example}

\newtheorem{rmks}[thm]{Remarks}

\title{Expansive algebraic actions of countable abelian groups} 
\author{Richard Miles}
\date{8 February, 2005}

\begin{document}

\maketitle

\begin{abstract}
This paper gives an algebraic characterization of expansive actions of 
countable abelian groups on compact abelian groups. This naturally extends the classification of 
expansive algebraic $\mathbb{Z}^d$-actions given by Schmidt using complex varieties. Also included is an application to 
a natural class of examples arising from unit subgroups of integral domains.
\end{abstract}


\thanks{The author is very grateful for the helpful advice and comments of Tom Ward.}

\section{Introduction}
Throughout, $G$ will be a countable abelian group, $X$ a compact metrizable abelian group and $\alpha$ an action of $G$ by automorphisms $\alpha_g$ of $X$. The dynamical system $(X,\alpha)$ is \emph{expansive} if and only if there is a neighbourhood $N$ of $0$ in $X$ such that
\begin{equation}\label{expansiveness_definition_formula}
\bigcap_{g\in G}\alpha_g(N) = \{0\}.
\end{equation}
Expansiveness is a property which has been studied both in its own right and in relation to other dynamical behaviours. Recently, Bhattacharya \cite{bhattacharya} has investigated expansiveness of arbitrary semi-groups on connected metrizable groups and has given a complete algebraic description of expansive actions on finite-dimensional connected abelian groups. The non-abelian case is covered by \cite[Theorem 2.4]{schmidt_book}. Also, using techniques applicable to Noetherian rings, Schmidt \cite{schmidt_expansiveness} has provided a classification of expansive $\mathbb{Z}^d$-actions on compact abelian groups in terms of complex varieties. The case for zero-dimensional groups has also been resolved \cite{kitchens_schmidt}, although some of the subtleties evident in \cite{schmidt_expansiveness} and \cite{bhattacharya} do not present themselves in this situation. Earlier studies of expansiveness for classical $\mathbb{Z}$-actions may be found, for example in \cite{eisenberg_vector_spaces}.

The expansive behaviour of algebraic actions has been exploited in relation to various other dynamical properties in settings such as \cite{ward_periodic}, \cite{schmidt_homoclinic} and \cite{einsiedler_schmidt}. Further investigations include the  introduction of a notion of expansive subdynamics  \cite{boyle_lind_esd}, developed in \cite{einsiedler_et_al_esd} and  \cite{einsiedler_ward}.

For an arbitrary $G$-action on a compact abelian group $X$, it is possible to realize the Pontryagin dual group $M=\widehat{X}$ as a countable module over the group ring $\mathbb{Z}[G]$, by defining \[
fa=\sum_{g\in G}c_g\widehat{\alpha}_g(a)
\]
where $a\in M$, $c_g\in\mathbb{Z}$ and $f=\sum_{g\in G}c_g g\in\mathbb{Z}[G]$ has $c_g=0$ for all but finitely many $g$. Conversely, any countable $\mathbb{Z}[G]$-module $M$ induces a natural action of $G$ on $X=\widehat{M}$ by setting $\alpha_g$ to be the automorphism of $X$ which is dual to multiplication by $g$ on $M$.
Roughly speaking, this duality turns topological properties into algebraic ones, an idea which plays a key role in the study of algebraic dynamical systems. Frequently, it will be necessary to consider quotients of the ring $\mathbb{Z}[G]$ and from now on the notation $\overline{G}$, $\overline{g}$ will be adopted to denote the respective images of $G$ and $g\in G$ under the natural quotient map. 

This paper gives an algebraic characterization of expansive actions of arbitrary countable abelian groups on compact abelian groups, as follows

\begin{mainresult}
Let $\alpha$ be an action of a countable abelian group $G$ by automorphisms of a compact abelian group $X$. Then $(X,\alpha)$ is expansive if and only if $M=\widehat{X}$ is a finitely generated $\mathbb{Z}[G]$-module and as $\mathfrak{a}$ runs through the annihilators of a set of generators for $M$, there is no ring homomorphism $\phi:\mathbb{Z}[G]/\mathfrak{a}\rightarrow\mathbb{C}$ for which $\phi(\overline{G})$ is a subgroup of the unit circle.
\end{mainresult}

The main difficulties involved in extending the approach used by Schmidt for $\mathbb{Z}^d$--actions (in particular \cite[Lemma 6.8]{schmidt_book}) relate to the fact that $\mathbb{Z}[G]$ may not be a Noetherian ring. This also means that the methods of decomposition employed in \cite{schmidt_book}, using associated prime ideals of Noetherian rings, are not applicable.

In Section \ref{valuations_section}, some algebraic consequences of the above theorem are investigated via a natural class of examples arising from unit subgroups of integral domains. In this setting, it is possible to exploit certain valuative maps to classify expansive behaviour.

\section{Preliminaries}
With $M=\widehat{X}$  realized as a module over $\mathbb{Z}[G]$, via duality, any submodule $N\subset M$ induces a quotient action of $G$ on a group of the form $X/Y$, where $Y$ is a closed $\alpha$-invariant subgroup of $X$. Similarly, a module of the form $M/N$ induces a $G$-action on a closed $\alpha$-invariant subgroup $Y$ of $X$ by restriction. 

The action of $G$ on $X$ can always be lifted to an action of $\Lambda=\bigoplus_\mathbb{N}\mathbb{Z}$, the group of eventually zero sequences of integers. By fixing a sequence of generators $(g_n)\in G^{\mathbb{N}}$, we induce a group homomorphism $\theta:\Lambda\rightarrow G$ defined by $\theta(m)=g^m$ where $m=(m_1,m_2,\dots)$ and $g^m=g_1^{m_1}g_2^{m_2}\dots$. Consequently, for each $m\in\Lambda$, setting 
\begin{equation}\label{lambda_action}
\alpha_m=\alpha_{\theta(m)}
\end{equation}
 induces an action of $\Lambda$ on $X$. It follows that for any $Y\subset X$ the intersections $\bigcap_{l\in\Lambda}\alpha_l(Y)$ and $\bigcap_{g\in G}\alpha_g(Y)$ agree. Thus, expansiveness of the $G$-action on $X$ can be studied in terms of the $\Lambda$-action described above. From now on the group ring $\mathbb{Z}[\Lambda]$ will be denoted by $R_\infty$. Dually, the viewpoint just described presents $\widehat{X}$ as an $R_\infty$-module and there is a ring homomorphism $\phi:R_\infty\rightarrow\mathbb{Z}[G]$ induced by sending $m\in\Lambda$ to $\theta(m)$. Furthermore, when $\mathfrak{a}$ is an ideal of $\mathbb{Z}[G]$, a $\mathbb{Z}[G]$-module of the form $\mathbb{Z}[G]/\mathfrak{a}$ can be considered as an $R_\infty$-module via the ring homomorphism $\pi\circ\phi$, where $\pi:\mathbb{Z}[G]\rightarrow\mathbb{Z}[G]/\mathfrak{a}$ denotes the natural map. In this case
\begin{equation}\label{canonical_module}
\widehat{X}\cong R_\infty/\mathfrak{b}
\end{equation}
where $\mathfrak{b}$ is the kernel of the map $\pi\circ\phi$. Multiplication by $m\in\Lambda$ in this ring is dual to the automorphism $\alpha_m$ given by (\ref{lambda_action}). In the sequel, quotient rings of the form just described will be particularly important. 

The isomorphism (\ref{canonical_module}) provides a canonical description of $X=\widehat{R_\infty/\mathfrak{b}}$ and the action of $\Lambda$ defined by (\ref{lambda_action}). By duality, we may identify $X$ with the subgroup of $\mathbb{T}^\Lambda$ consisting of elements $x=(x_l)$ satisfying $\langle x,f\rangle=1$ for all $f\in\mathfrak{b}$. Let
\begin{equation}\label{group_ring_element}
f = \sum_{l\in\Lambda}c_l(f)l
\end{equation}
where $c_l(f)=0$ for all but finitely many $l\in\Lambda$. Then
\[
\langle x,f\rangle = \exp 2\pi i{\sum_{l\in\Lambda}c_l(f)x_l}.
\] 
Hence, $X$ consists of elements $x=(x_l)\in\mathbb{T}^\Lambda$ satisfying
\begin{equation}\label{annihilator_condition}
\sum_{l\in\Lambda}c_l(f)x_{l+m} = 0 \mod 1
\end{equation}
for all $f\in\mathfrak{b}$ and $m\in\Lambda$. Using this description of $X$, $\alpha_m$ can be interpreted as the restriction of the shift map $(x_l)\mapsto(x_{l+m})$ to $X$. This defines the \emph{natural action} of $\Lambda$ on $X$.

\section{Algebraic Criteria}\label{criteria_section}

Let $B=\ell^\infty(\Lambda)$ denote the Banach space of all bounded complex-valued functions $z=(z_l)$ on $\Lambda$ with the supremum norm. For each $m\in\Lambda$, the shift map $\sigma_m:B\rightarrow B$ given by
\begin{equation}\label{shift_map_definition}
\sigma_m\left((z_l)\right) = (z_{l+m})
\end{equation}
is an isometry of $B$. Any $f\in R_\infty$ defines a linear operator $T_f$ on $B$ as follows
\begin{equation}\label{lin_op_of_f}
T_f(z) = \sum_{m\in\Lambda}c_m(f)\sigma_m(z)
\end{equation}
where $c_m(f)$ is given by (\ref{group_ring_element}). 

\begin{lem}\label{enough_points_lemma}
Let $X=\widehat{R_\infty/\mathfrak{a}}$ and suppose that $\alpha$ is the natural action of $\Lambda$ on $X$. If $(X,\alpha)$ is non-expansive then for any finite subset $F$ of $\mathfrak{a}$, there is a non-trivial $z\in\ell^\infty(\Lambda)$ such that for all $f\in F$, $T_f(z) = 0$, where $T_f$ is given by (\ref{lin_op_of_f}).
\end{lem}

\begin{proof}
First identify $X$ with the shift-invariant subgroup of $\mathbb{T}^{\Lambda}$ given by (\ref{annihilator_condition}) and $\alpha$ with the restriction of the shift. Let $K=\max \left\{\sum_{l\in\Lambda}|c_l(f)|:f\in F\right\}$ and $\varepsilon=(10K)^{-1}$. Set $N=\left\{x=(x_l)\in X:\|x_0\|<\varepsilon\right\}$ where $\|\cdot\|$ is given by $\|t\|=\min\left\{|a-t|:a\in\mathbb{Z}\right\}$. Since $(X,\alpha)$ is non-expansive, $\bigcap_{m\in\Lambda}\alpha_m(N)$ contains a non-zero $x\in X$. Moreover, the identification of $\alpha_m$ with the shift means 
$\|x_l\|<\varepsilon$ for all $l\in\Lambda$. Choose $z=(z_l)\in \ell^\infty(\Lambda)$ such that $z_l=x_l$,  $l\in\Lambda$. Since $x$ satisfies (\ref{annihilator_condition}), it follows that for all $m\in\Lambda$ and $f\in F$, $\sum_{l\in\Lambda}c_l(f)z_{l+m} \in \mathbb{Z}$. Hence, by the choice of $\varepsilon$ 
\[
\left|\sum_{l\in\Lambda}c_l(f)z_{l+m}\right| < 1
\]
which means $\sum_{l\in\Lambda}c_l(f)z_{l+m}=0$ for all $m\in\Lambda$ and $f\in F$. Equivalently, $T_f(z)=0$ for all $f\in F$.
\end{proof}

\begin{lem}\label{cyclic_module_lemma}
Let $X=\widehat{R_\infty/\mathfrak{a}}$ and suppose that $\alpha$ is the natural action of $\Lambda$ on $X$. If $(X, \alpha)$ is non-expansive then there is a ring homomorphism $\phi:R_\infty/\mathfrak{a}\rightarrow\mathbb{C}$ for which $\phi(\overline{\Lambda})$ is a subgroup of the unit circle.
\end{lem}

\begin{proof}
Choose an enumeration $\mathfrak{a}=\{f_1,f_2,\dots\}$. Let $B=\ell^\infty(\Lambda)$ and for each $k\in\mathbb{N}$, set $F_k = \{f_1,\dots,f_k\}$ and
\[
B_k=\left\{z\in B:T_f(z) = 0\textrm{ for every }f\in F_k\right\}.
\]
Then
\[
B_1\supset B_2\supset B_3 \supset \dots
\]
is a chain of closed linear subspaces of $B$. Furthermore, each $B_k$ is non-trivial by Lemma \ref{enough_points_lemma}. Let $\mathcal{A}$ denote the Banach algebra of bounded linear operators on $B$ and $\mathcal{B}$ the sub-algebra generated by $\{\sigma_m:m\in\Lambda\}$ where $\sigma_m\in\mathcal{A}$ is given by (\ref{shift_map_definition}). Similarly, for each $k\in\mathbb{N}$, let $\mathcal{B}_k$ denote the algebra of bounded linear operators on $B_k$ generated by $\{\sigma_m|_{B_k}:m\in\Lambda\}$. Restricting the shift maps in this way induces ring epimorphisms $\mathcal{B}\rightarrow\mathcal{B}_k$ with corresponding 
kernels $\mathfrak{b}_1\subset\mathfrak{b}_2\subset\mathfrak{b}_3\subset\dots$. Note that $T_f\in\mathfrak{b}_k$ for all $f\in F_k$. Let $\mathfrak{b}\subset\mathcal{B}$ be the ideal generated by $\{T_f:f\in\mathfrak{a}\}$. If $1\in\mathfrak{b}$ then for some $k$, $1\in\mathfrak{b}_k$  which implies $\mathcal{B}_k=\{0\}$ and this contradicts the definition of $B_k$. Therefore, $\mathfrak{b}$ is a proper ideal of $\mathcal{B}$. Hence there is a maximal ideal $\mathfrak{m}\subset\mathcal{B}$ containing $\mathfrak{b}$. 

Let $\mathcal{M}(\mathcal{B})$ denote the maximal ideal space of $\mathcal{B}$, $\mathcal{C}(\mathcal{M}(\mathcal{B}),\mathbb{C})$ the space of continuous complex-valued functions on $\mathcal{M}(\mathcal{B})$ and  $\hat{}:\mathcal{B}\rightarrow\mathcal{C}(\mathcal{M}(\mathcal{B}),\mathbb{C})$ the Gelfand transform. The choice of $\mathfrak{m}$ means that 
\begin{equation}\label{gelfand_transform_is_zero}
\hat{T}_f(\mathfrak{m}) = \sum_{m\in\Lambda}c_m(f)\hat{\sigma}_m(\mathfrak{m}) = 0
\end{equation}
for all $f\in\mathfrak{a}$. For each $n\in\mathbb{N}$, let $\lambda_n=\hat{\sigma}_{e(n)}(\mathfrak{m})$ where $e(n)$ is the $n$-th unit vector in $\Lambda$ and set $w=(\lambda_n)\in\mathbb{C}^{\mathbb{N}}$. For every  $f\in\mathfrak{a}$, it follows that
\[
\sum_{m\in\Lambda}c_m(f)w^m = \sum_{m\in\Lambda}c_m(f)\hat{\sigma}_m(\mathfrak{m}) 
\]
where $w^m=\lambda_1^{m_1}\lambda_2^{m_2}\dots$, $m=(m_1,m_2,\dots)$. Moreover, (\ref{gelfand_transform_is_zero}) implies that $w$ induces a well defined evaluation homomorphism $\phi_w:R_\infty/\mathfrak{a}\rightarrow\mathbb{C}$ given by substituting $w^l$ for $l$ in the expression (\ref{group_ring_element}). Finally, since each $\sigma_m$ is an isometry of $B$, $\|\sigma_m\|=1$ and the Gelfand transform ensures that $|w^m|=1$ for all $m\in\Lambda$. It follows that $\phi_w(\overline{\Lambda})$ is a subgroup of the unit circle.
\end{proof}

\begin{thm}\label{main_theorem}
Let $\alpha$ be an action of a countable abelian group $G$ by automorphisms of a compact abelian group $X$. Then $(X,\alpha)$ is expansive if and only if $M=\widehat{X}$ is a finitely generated $\mathbb{Z}[G]$-module and as $\mathfrak{a}$ runs through the annihilators of a set of generators for $M$, there is no ring homomorphism $\phi:\mathbb{Z}[G]/\mathfrak{a}\rightarrow\mathbb{C}$ for which $\phi(\overline{G})$ is a subgroup of the unit circle.
\end{thm}
\begin{proof}
If $M$ is not finitely generated then \cite[Proposition 4.3]{bhattacharya} shows that $(X,\alpha)$ cannot be expansive. Hence assume that $M$ is finitely generated and there is a generator $a$ of $M$ whose annihilator $\mathfrak{a}\subset\mathbb{Z}[G]$ is such that there exists
a ring homomorphism $\phi:\mathbb{Z}[G]/\mathfrak{a}\rightarrow\mathbb{C}$ with $\phi(\overline{G})$ a subgroup of the unit circle.
Note that $\phi$ is necessarily induced by a ring homomorphism $\Phi:\mathbb{Z}[G]\rightarrow\mathbb{C}$ with $\Phi(G)=\phi(\overline{G})$. Furthermore, $\mathbb{C}$ can be considered as a $\mathbb{Z}[G]$-module via $\Phi$ which can be used to construct a module homomorphism $\theta:M\rightarrow\mathbb{C}$ as follows. 

Suppose $N$ is a submodule of $M$ and there is a module homomorphism $\theta_N:N\rightarrow\mathbb{C}$. A submodule of this form exists since we can set $N=\mathbb{Z}[G]a$ and $\theta_N(fa)=\Phi(f)$, $f\in\mathbb{Z}[G]$. Let $\mathfrak{p}=\ker \Phi$ and observe that $\theta_N$ is well defined, since by the definition of $\Phi$, $\mathfrak{a}$ must be contained in $\mathfrak{p}$.
Let $b\in M\setminus N$ be a generator of $M$ and set $L=N+\mathbb{Z}[G]b$. To see that there is also a module homomorphism $\theta_L:L\rightarrow\mathbb{C}$, consider the following two cases. If $f\in\mathbb{Z}[G]$ and $fb\in N$ implies $f\in\mathfrak{p}$ then define $\theta_L:L\rightarrow\mathbb{C}$  by
\[
\theta_L(c+fb)=\Phi(f)
\] 
where $c\in N$, $f\in\mathbb{Z}[G]$. Alternatively, if there exists $h\in\mathbb{Z}[G]\setminus\mathfrak{p}$ with $hb\in N$, set $\lambda=\theta_N(hb)/\Phi(h)$ and define $\theta_L:L\rightarrow\mathbb{C}$  by
\[
\theta_L(c+fb)=\theta_N(c)+\Phi(f)\lambda
\] 
where $c\in N$ and $f\in\mathbb{Z}[G]$. In both cases it is routine to check that $\theta_L$ is well defined. It now follows by induction on the generators of $M$ that there must be a module homomorphism $\theta=\theta_M$ from $M$ to $\mathbb{C}$. 

The closure of $\theta(M)$ in $\mathbb{C}$, $W$ say, is a submodule of $\mathbb{C}$ and $\Phi$ induces an action of $G$ on $W$ via multiplication by $\Phi(g)$, $g\in G$. Denote the corresponding dual action on $\widehat{W}$ by $\beta$. Since $\Phi(G)$ is a subgroup of the unit circle, $(\widehat{W}, \beta)$ is non-expansive. Also $\theta(M)$ is dense in $W$ and this means that there is an induced equicontinuous injection $\widehat{\theta}:\widehat{W}\rightarrow X$. Moreover, since $\theta$ is a module homomorphism, by duality there is a conjugacy $\alpha_g\circ\widehat{\theta}=\widehat{\theta}\circ\beta_g$, $g\in G$. 
This means the restriction of $\alpha$ to $\widehat{\theta}(\widehat{W})$ is non-expansive. Thus $(X,\alpha)$ must be non-expansive.

Conversely, assume that $M$ is finitely generated by $k$ elements with corresponding annihilators $\mathfrak{a}_1,\dots,\mathfrak{a}_k$. For $1\leq j \leq k$, set $M_j=\mathbb{Z}[G]/\mathfrak{a}_j$ and suppose for each $j$, there is no ring homomorphism $\phi:M_j\rightarrow\mathbb{C}$ for which $\phi(\overline{G})$ is a subgroup of the unit circle. Employing the realization of cyclic modules provided by (\ref{canonical_module}), Lemma \ref{cyclic_module_lemma} shows that the induced action of $G$ on each $\widehat{M_j}$ is expansive. Thus the product action on $Y=\widehat{M_1}\times\dots\times\widehat{M_k}$ is expansive. The corresponding dual $\mathbb{Z}[G]$-module is $N=M_1\oplus\dots\oplus M_k$ and there is a surjective module homomorphism of $N$ onto $M$. Consequently, by duality $X$ may be regarded as a closed subgroup of $Y$ and $\alpha$ may be identified with the corresponding restriction of the product action. Thus $(X,\alpha)$ is expansive.  
\end{proof}

The following example is a simple (but non-trivial) illustration of a non-expansive $\Lambda$-action arising from a cyclic module over the group ring $R_\infty$. 

\begin{eg}\label{non_exp_cyclic_module_eg}
For each $n\in\mathbb{N}$, let $e(n)$ be the $n$-th unit vector in $\Lambda$.  Consider the proper ideal $\mathfrak{a}\subset R_\infty$ generated by elements of the form
\[
f_n = e(n+1)-(n+1)e(1)+n 
\]
where $n\in\mathbb{N}$. Let $M=R_\infty/\mathfrak{a}$. The corresponding natural action $\alpha$ of $\Lambda$ on $X=\widehat{M}$ cannot be expansive. To see this, define a ring homomorphism $\Phi:R_\infty\rightarrow\mathbb{C}$ by setting $e(n)\mapsto 1$, $n\in\mathbb{N}$. Note that $\Phi(\Lambda)=\{1\}$. Since $\Phi(f_n)=0$ for every $n\in\mathbb{N}$, there is an induced ring homomorphism $\phi:M\rightarrow\mathbb{C}$ with $\phi(\overline{\Lambda})=\{1\}$. Hence by Theorem \ref{main_theorem}, the dynamical system $(X,\alpha)$ is non-expansive. For a more direct perspective, consider $X$ as a subgroup of $\mathbb{T}^{\Lambda}$ consisting of elements $(x_l)$ satisfying (\ref{annihilator_condition}) for every $f=f_n$ and $m\in\Lambda$. Given any $t\in\mathbb{T}$, setting $x_l=t$ for all $l\in\Lambda$ determines an element $(x_l)\in \mathbb{T}^{\Lambda}$ satisfying (\ref{annihilator_condition}). Such shift invariant points in $X$ readily demonstrate the non-expansiveness of this dynamical system.
\end{eg}

A straightforward consequence of the algebraic characterization of expansiveness for $\mathbb{Z}^d$-actions is that a factor of an expansive $\mathbb{Z}^d$-action is again expansive \cite[Corollary 6.15]{schmidt_book}. The following example shows that this need not be the case in the more general setting considered here. The author would like to thank the referee for pointing out that such an example might exist.

\begin{eg}
Let $e(n)$, $n\in\mathbb{N}$ be as in Example \ref{non_exp_cyclic_module_eg}. Set $\mathfrak{a}=(e(1)-2)\subset R_\infty$. Then any ring homomorphism from the cyclic $R_\infty$-module $M=R_\infty/\mathfrak{a}$ to $\mathbb{C}$ cannot send $\overline{\Lambda}$ to a subgroup of the unit circle. Hence the natural $\Lambda$-action $\alpha$ on $X=\widehat{M}$ is expansive. With $M$ considered as a ring, any ideal $\mathfrak{b}\subset M$ is a submodule of $M$. For example, take the proper ideal $\mathfrak{b}$ generated by elements of the form $\overline{e(n)}+\overline{e(n+1)}$, where $n\geq 2$ is even. Now $\mathfrak{b}$ is not finitely generated as an $R_\infty$-module and so the natural $\Lambda$-action $\beta$ on $Y=\widehat{\mathfrak{b}}$ is not expansive. Moreover, by duality the non-expansive system $(Y,\beta)$ is a factor of the expansive system $(X,\alpha)$.   
\end{eg}

\begin{rmks}
For the generators $f_n$ of $\mathfrak{a}$ in Example \ref{non_exp_cyclic_module_eg}, $\sum_{l\in\Lambda}|c_l(f_n)|$ can be arbitrarily large. It is ideals of this nature, that arise in the non-Noetherian ring $R_\infty$, which make the adaptation of the proof of \cite[Lemma 6.8]{schmidt_book} more difficult in our setting.

A ring homomorphism from a $\mathbb{Z}[G]$-module of the form $M_G=\mathbb{Z}[G]/\mathfrak{a}$ to $\mathbb{C}$ is determined by its image of $G$ in $\mathbb{C}$. The collection of all such homomorphisms may be identified with the $\mathbb{C}$-valued points of the affine scheme $\textrm{Spec }M_G$ (see for example \cite[Section 1.4]{eisenbud_harris}). If $G\cong\mathbb{Z}^d$ and $\mathfrak{a}$ does not contain a non-zero constant, then the the set of $\mathbb{C}$-valued points is a complex variety. The classification of expansive $\mathbb{Z}^d$-actions due to Schmidt \cite{schmidt_expansiveness} uses this fact. 

It is well known that an expansive algebraic $\mathbb{Z}^d$-action $\alpha$ on a compact abelian group $X$ must satisfy the descending chain condition on closed $\alpha$-invariant subgroups of $X$. That is every chain 
\[
X\supset X_1\supset X_2\supset\dots
\]
of closed $\alpha$-invariant subgroups, eventually becomes stationary. For acting groups like $\Lambda$, this is not necessarily the case, even if $M=\widehat{X}$ is a cyclic $R_\infty$-module. As a simple example, consider a sequence of distinct rational primes $\{p_n\}$ and a corresponding sequence of ideals $\{\mathfrak{a}_n\}$ in $R_\infty$, where  $\mathfrak{a}_n$ is generated by $e(1)-p_1,\dots,e(n)-p_n$, with $e(n)$ as in Example \ref{non_exp_cyclic_module_eg}. Setting $M_n=R_\infty/\mathfrak{a}_n$, we obtain a sequence of $R_\infty$-modules and natural projections
\[
M_1 \twoheadrightarrow M_2 \twoheadrightarrow \dots  
\]
Dually, this gives a sequence of shift-invariant subgroups of $\mathbb{T}^{\Lambda}$
\begin{equation}\label{descending_chain_example}
X_1 \leftarrow X_2 \leftarrow \dots
\end{equation}
where $X_n=\widehat{M_n}$ and each arrow represents the natural inclusion of $X_{n+1}$ in $X_n$. The construction of the $\mathfrak{a}_n$ means that these inclusions are proper and the chain (\ref{descending_chain_example}) does not become stationary. Moreover, using Theorem \ref{main_theorem}, it is routine to check that the $\Lambda$-shift on $X_1$ is expansive.
\end{rmks}

\section{Valuations and Expansiveness}\label{valuations_section}
In this section, the results of Section \ref{criteria_section} are applied to a natural class of dynamical systems that allow a characterization of expansiveness using valuations. Suppose that $M$ is an integral domain and $G$ is a subgroup of the unit group of $M$. Then $M$ is naturally a $\mathbb{Z}[G]$-algebra and the map dual to multiplication by $g\in G$ on $M$ is an automorphism of $X=\widehat{M}$. This induces a $G$-action $\alpha$ by automorphisms of $X$. Furthermore, the corresponding representation of $\widehat{X}$ as a $\mathbb{Z}[G]$-module agrees with the natural $\mathbb{Z}[G]$-algebra structure of $M$. From now on, attention is restricted to finitely generated $\mathbb{Z}[G]$-algebras $M$. If $M$ is not finitely generated as an algebra over $\mathbb{Z}[G]$, then it cannot be finitely generated as a module and Theorem \ref{main_theorem} implies that the corresponding $G$-action is non-expansive. Understanding both the structure of the group $G$ and the structure of $M$ as a $\mathbb{Z}[G]$-algebra is at the heart of determining whether or not the dynamical system $(X,\alpha)$ is expansive. 

In what follows, two types of valuative function will be important. Firstly, for any domain $R$ and homomorphism $\phi$ of $R$ into the complex numbers, there is a corresponding \emph{logarithmic map} $v_\phi:R\rightarrow\mathbb{R}\cup\{\infty\}$ given by 
\[
v_\phi(a)=\log|\phi(a)|
\]
where $a\in R$. Secondly, we will consider genuine valuations $v:R\rightarrow H\cup\{\infty\}$, where $H$ is an ordered group and $v$ corresponds to a valuation ring in the field of fractions of $R$. The relationship between logarithmic maps and valuations is discussed in \cite{bergman}.

We now return to our domain $M$ and subgroup $G$ of the unit group. Since $M$ is a domain, when it is considered as a $\mathbb{Z}[G]$-module, there is a unique annihilator $\mathfrak{p}\subset\mathbb{Z}[G]$ corresponding to every element of $M$. This is simply the kernel of the natural homomorphism from $\mathbb{Z}[G]$ to $M$ and $\mathbb{Z}[G]/\mathfrak{p}$ is isomorphic to the smallest subring of $M$ containing $G$, which from now on will be denoted by $M_G$. 
Denote the collection of all logarithmic maps $v:M_G\rightarrow\mathbb{R}\cup\{\infty\}$ arising from homomorphisms of $M_G$ into $\mathbb{C}$ by $\mathcal{V}_\infty(M_G)$.
Now consider the collection of all valuation rings of the field of fractions $k$ of $M$. Denote the Zariskii Space (or Zariskii Manifold) of valuation rings of $k$ which contain a subring $R\subset k$ by $\textrm{Zar}(k,R)$. The space $\textrm{Zar}(k,R)$ is a multifaceted object upon which a natural topology may be defined \cite[Chapter 10]{matsumura}. Set $R=\mathbb{Z}$ or $R=\mathbb{F}_p$ according to the characteristic of $k$ and let $\mathcal{Z}=\textrm{Zar}(k,R)\setminus\textrm{Zar}(k,M)$. Denote the corresponding collection of valuations by $\mathcal{V}_\mathcal{Z}(M)$. 

\begin{prop}\label{main_valuation_theorem}
Let $M$ be an integral domain and $G$ a subgroup of its unit group such that $M$ is a finitely generated $\mathbb{Z}[G]$-algebra. Then the induced action of $G$ on $X=\widehat{M}$ is expansive if and only if $v|_G$ is non-trivial for all $v\in\mathcal{V}_\infty(M_G)\cup\mathcal{V}_\mathcal{Z}(M)$. 
\end{prop}

\begin{proof}
Theorem \ref{main_theorem} together with the above discussion shows that the induced action of $G$ on $M$ is expansive if and only if $M$ is finitely generated as a $\mathbb{Z}[G]$-module and there is no ring homomorphism $\phi:M_G\rightarrow\mathbb{C}$ for which $\phi(G)$ is a subgroup of the unit circle. This latter condition is equivalent to $v(G)\neq\{0\}$ for all $v\in\mathcal{V}_\infty(M_G)$. Hence it remains to show that $M$ being finitely generated as a $\mathbb{Z}[G]$-module is equivalent to $v(G)\neq\{0\}$ for all $v\in\mathcal{V}_\mathcal{Z}(M)$. Since $M$ is finitely generated as an algebra, \cite[Theorem 9.1]{matsumura} shows that $M$ is finitely generated as a $\mathbb{Z}[G]$-module if and only if every element of $M$ is integral over $M_G$. That is, the integral closures of $M$ and $M_G$ coincide in $k$, the field of fractions of $M$. By \cite[Theorem 10.4]{matsumura} this is equivalent to $\mathcal{Z}\cap\textrm{Zar}(k,M_G)=\varnothing$. This happens precisely when $v(G)\neq\{0\}$ for all $v\in\mathcal{V}_\mathcal{Z}(M)$.
\end{proof}

\begin{eg}\label{simple_field_example}
Let $M$ be a countable field, $X=\widehat{M}$ and $G$ a non-trivial subgroup of $M^{\times}$. Let $\alpha$ denote the induced action of $G$ on $X$ and 
consider $\mathcal{V}_\mathcal{Z}(M)$. This contains all non-trivial valuations on the field $M$. We need only find a single valuation $v$ such that $v|_G$ is trivial for the dynamical system $(X,\alpha)$ to be non-expansive. For example, if $G$ is a subgroup of the unit group of a valuation ring in $M$ with corresponding valuation $v$, then $v(G)=\{0\}$. 
\end{eg}

\begin{eg}
Let $M$ and $X$ be as in Example \ref{simple_field_example}. If $M$ has positive characteristic then $\mathcal{V}_\infty(M_G)=\varnothing$. If $M$ has zero characteristic and $G$ contains an integer $n\neq 1$, then any $v\in\mathcal{V}_\infty(M_G)$ has $v(n)=\log|n|$. Hence in either situation, expansiveness of $(X,\alpha)$ depends only on $\mathcal{V}_\mathcal{Z}(M)$. 
If $G=M^{\times}$ then $M_G=M$ and so $(X,\alpha)$ is expansive. In addition, Example \ref{simple_field_example} gives a rich supply of sub-actions which are non-expansive. That is subgroups $G\subset M^{\times}$ for which the natural action of $G$ on $X$ is non-expansive. 
On the other hand, if a subfield $L$ is chosen such that $L\subset M$ is a finite extension and $G=L^{\times}$, 
then $(X,\alpha)$ is expansive. To see this in terms of $\mathcal{V}_\mathcal{Z}(M)$, note that
any valuation $v$ which is trivial on $G$ must also be trivial on $L$ and hence $M$.  However, $\mathcal{V}_\mathcal{Z}(M)$ is assumed not to contain the trivial valuation.  
\end{eg}

\begin{eg}
Consider the function field $k=\mathbb{Q}(t)$ and a subgroup $H\subset k^{\times}$ generated by (at least 2) distinct irreducible polynomials of degree 1 in $\mathbb{Z}[t]$. For simplicity, assume $t\in H$. Let $M$ be the smallest subring of $k$ containing $H$ and suppose $G$ is a subgroup of $H$. If $G\neq H$ then the induced action of $G$ on $X=\widehat{M}$ cannot be expansive. To see this, simply take an irreducible $f\in H\setminus G$ and localise $\mathbb{Z}[t]$ at the prime ideal generated by $f$. This gives a discrete valuation ring in $k$ whose maximal ideal contains $f$. Furthermore, there is a corresponding valuation $v\in\mathcal{V}_\mathcal{Z}(M)$ with $v|_G$ trivial. 

Now suppose $G=H$. Then $M_G=M$ and it is only necessary to consider $\mathcal{V}_\infty(M_G)$ to check expansiveness. Suppose there is a logarithmic map $v$ such that $v(G)=\{0\}$. Then $v(t)=0$ and $v(at+b)=0$ for integers $a\neq b$. Treating these equalities simultaneously, it follows that this can happen if and only if
\[
a^4-2(b^2+1)a^2+(b^2-1)^2\leq 0
\]
Since $a$ and $b$ are integers, the only possibilities are $a=\pm(b+1)$, $a=\pm(b-1)$ or $a=\pm{b}$. In general, if the set of generators for $G$ is confined to irreducibles arising from finitely many distinct $a$, $G$ can of course still be infinitely generated. However, in this case, if $G$ acts non-expansively, the restricted choices of $b$ mean that in fact $G$ must be be finitely generated. If infinitely many distinct $a$ are permitted, then appropriate choices of $b$ allow us to construct both expansive and non-expansive $G$-actions, even when $M$ is not finitely generated as a ring over $\mathbb{Z}$.
\end{eg}

\begin{rmks}
In general $\mathcal{V}_\mathcal{Z}(M)$ can be complicated \cite[Section 1.7]{cohns_book}, except in special circumstances, say where $M$ is a ring of algebraic integers. Therefore, alternative approaches can be equally useful for determining the expansiveness or otherwise of the dynamical systems described here. For example, a direct application of the integrality criterion contained in the proof of Proposition \ref{main_valuation_theorem} can be useful, depending on how the domain $M$ is presented.

For $G$ finitely generated, the sets $\mathcal{V}_\infty(M_G)$ and $\mathcal{V}_\mathcal{Z}(M)$ were investigated by the author in \cite{miles_thesis} to study the expansive subdynamics of algebraic $\mathbb{Z}^d$-actions, developed in \cite{einsiedler_et_al_esd}. 
In this case, $\mathcal{V}_\infty(M_G)$ produces the logarithmic image of a complex variety (or amoeba) $A\subset\mathbb{R}^d$. To complete the picture, $\mathcal{V}_\mathcal{Z}(M)$ can be replaced by a collection of discrete valuations yielding a `non-archimedean' companion to $A$ \cite[Theorem 4.3.10]{miles_thesis}. 
There are of course other situations where $\mathcal{V}_\mathcal{Z}(M)$ may profitably be replaced by a collection of discrete valuations. In general, it makes sense to ask when $\mathcal{V}_\mathcal{Z}(M)$ can be restricted in some way, without affecting the statement of Proposition \ref{main_valuation_theorem}. 
In a slightly different direction, Einsiedler, Kapranov and Lind \cite{einsiedler_kapranov_lind} have explained important links between non-archimedean amoebas and the set $\mathcal{V}_\mathcal{Z}(M)$. Their approach elucidates the role played by valuations in the study of expansive subdynamics of algebraic $\mathbb{Z}^d$-actions. 

Finally, it should be noted that a consistent notion of expansive subdynamics for acting groups like $\Lambda$ is yet to be established.

\end{rmks}

Richard Miles, University of East Anglia, Norwich, NR4 7TJ, UK\\
\indent\emph{E-mail address}: \texttt{r.miles@uea.ac.uk}

\end{document}